\documentclass[titlepage,11pt]{article}
\usepackage{latexsym}
\usepackage{amsfonts}
\usepackage{amsmath}
\usepackage{mathtools}

\usepackage{tikz}
\usetikzlibrary {shapes.geometric}
\newcommand\blackslug{\hbox{\hskip 1pt \vrule width 4pt height 8pt depth 1.5pt
        \hskip 1pt}}
\newcommand\bbox{\hfill \quad \blackslug \bigbreak}

\def\CC{\hbox{-}\cdots\hbox{-}}
\def\LL{,\ldots,}

\newcommand{\cupcup}{\cup \cdots\cup}

\def\dist{\operatorname{dist}}

\title{Asymptotic structure. IV. A counterexample to the weak coarse Menger conjecture}
\author{
Tung Nguyen\thanks{Supported by AFOSR grants
A9550-19-1-0187 and FA9550-22-1-0234, and by NSF grants  DMS-1800053 and DMS-2154169.}\\
Princeton University,\\ Princeton, NJ 08544, USA
\and
Alex Scott\thanks{Supported by EPSRC grant EP/X013642/1}\\
University of Oxford, \\
Oxford, UK
\and
Paul Seymour\thanks{Supported by AFOSR grants
A9550-19-1-0187 and FA9550-22-1-0234, and by NSF grants  DMS-1800053 and DMS-2154169.}\\
Princeton University,\\ Princeton, NJ 08544, USA}

\date{July 20, 2025; revised \today}

\newtheorem{thm}{}[section]

\newcommand{\Proof}{\noindent{\bf Proof.}\ \ }

\begin{document}
\maketitle
\begin{abstract}
Coarse graph theory concerns finding ``coarse'' analogues of graph theory theorems, replacing disjointness with being far apart. One
of the most interesting open questions is to find a coarse analogue of Menger's theorem, which characterizes when there are $k$ 
vertex-disjoint paths between two given sets $S,T$ of vertices of a graph. We showed in an earlier paper that the most natural such analogue is false,
but a weaker statement remained as a popular open question. Here we show that the weaker statement is also false. 

More exactly, suppose that $S,T$ are sets of vertices of a graph $G$, and there do not exist $k$ paths between $S,T$, pairwise at 
distance at least $c$. To make an analogue of Menger's theorem, one would like to prove that there must be a small set 
$X\subseteq V(G)$ 
such that every $S-T$ path of $G$ passes close to a member of $X$: but how small and how close? In view of Menger's theorem,
one would hope for $|X|<k$ and ``close'' some function of $k,c$ (and indeed, this was conjectured by 
Georgakopoulos and Papasoglu, and  independently, by
Albrechtsen, Huynh, Jacobs, Knappe and Wollan); but we showed that this is false, even if $c=3$ and $k=3$. 

Here we upgrade the counterexample:
we show that, even if $c=k=3$, no pair of constants (for ``small'' and ``close'') work. For all $\ell, m$, there is a graph $G$ and $S,T\subseteq V(G)$,  such that
there do not exist three $S-T$ paths pairwise with distance at least three, and yet there is no $X$ with $|X|\le m$ 
such that
every $S-T$ path passes within distance at most $\ell$ of $X$. 

\end{abstract}

\section{Introduction}

The ``disjoint paths problem'' asks when there is a set of $k$ vertex-disjoint paths between sets $S,T$ of vertices of a graph $G$; 
and it is answered by a theorem of K. Menger from 1927~\cite{menger}, that such paths exist if and only if there is no subset $X\subseteq V(G)$ 
of size $<k$ such that every $S-T$ path has a vertex in $X$. But what if we want the paths to be at least a certain distance from one 
another?\footnote{If $X$ and $Y$ are vertices, or sets of vertices, or subgraphs, of a graph $G$, then $\dist_G(X,Y)$ denotes the distance between 
$X,Y$, that is, the number of edges in the shortest path of $G$ with one end in $X$ and the other in $Y$.}
This question is motivated both by the developing area of ``coarse graph theory'', which is concerned with the large-scale geometric 
structure of graphs (see Georgakopoulos and Papasoglu~\cite{agelos}), and by the algorithmic question of deciding whether such paths exist (see
Bienstock~\cite{bienstock}, Kawarabayashi and Kobayashi~\cite{kk}, and Balig\'acs and MacManus~\cite{bm}).

A coarse analogue of Menger's theorem was conjectured by
Albrechtsen, Huynh, Jacobs, Knappe and Wollan~\cite{wollan}, and independently by Georgakopoulos and Papasoglu~\cite{agelos}:
\begin{thm}\label{conj}
{\bf False conjecture:} For all integers $k,c\ge 1$ there exists $\ell>0$ with the following property.
Let $G$ be a graph and let $S,T\subseteq V(G)$.  Then either 
\begin{itemize}
\item there are $k$ paths between $S,T$, pairwise at distance at least $c$; or 
\item there is a set $X\subseteq V(G)$
with $|X|\le k-1$ such that every path between $S,T$ contains a vertex with distance at most $\ell$ from some member of $X$.
\end{itemize}
\end{thm}
Both sets of authors proved the result for $k=2$, but we showed in~\cite{counterex} that \ref{conj} is false for all $k\ge 3$,
even if $c=3$ and $G$ has maximum degree three. 
The case $c=3$ is of special interest, because it is easy to see that if the result is true when $c=3$ (for some value of $k$) then it 
is true for all 
$c\ge 3$ and the same value of $k$ (apply the result when $c=3$ to the $c$th power of $G$). And indeed, the conjecture remains open 
when $c=2$, and we have nothing to say about that case in this paper.

Since the most natural extension of Menger's theorem is false, we fall back onto what seems the next most natural, 
the following weaker statement:
\begin{thm}\label{weakcoarseconj}
{\bf False conjecture:} For all integers $k,c\ge 1$ there exist $m,\ell>0$ with the following property.
Let $G$ be a graph and let $S,T\subseteq V(G)$.  Then either
\begin{itemize}
\item there are $k$ paths between $S,T$, pairwise at distance at least $c$; or
\item there is a set $X\subseteq V(G)$
with $|X|\le m$ such that every path between $S,T$ contains a vertex with distance at most $\ell$ from some member of $X$.
\end{itemize}
\end{thm}
We proposed this in~\cite{counterex}, but 
now we will show that this too is false, even if $c=k=3$. More exactly, we will show:

\begin{thm}\label{mainthm}
For all integers $\ell,m\ge 1$, there is a graph $G$ and subsets $S,T\subseteq V(G)$ such that:
\begin{itemize}
\item there do not exist three paths between $S,T$ that pairwise have distance at least three; and
\item for every set $X\subseteq V(G)$ with $|X|\le m$, there is a path $P$ between $S,T$ such that $\dist_G(P,X)>\ell$.
\end{itemize}
\end{thm}

Our counterexample has some vertices of large degree, but it can easily be modified into a counterexample with only one vertex of degree more than three, while
keeping $c=k=3$. We will explain this in section 3.


\section{The counterexample}

For all integers $\ell,m\ge 1$, we will give a construction for a quadruple $(G,r,S,T)$ called (in this paper)  
an {\em $(\ell,m)$-block}, where $G$ is a graph, $r\in V(G)$, and $S,T$ are disjoint subsets of $V(G)\setminus \{r\}$, both of size $m$.
Later we will prove (by induction on $m$) that it has the properties that:
\begin{itemize}
\item $\dist_G(u,v)> 2\ell$ for every two vertices $u,v\in \{r\}\cup S\cup T$; 
\item for every choice of $X\subseteq V(G)$ with $|X|<m$, there is a path $P$ between $S,T$ such that $\dist_G(P,X \cup \{r\})>\ell$; and
\item for every two vertex-disjoint $S-T$ paths $P,Q$, either one of $P,Q$ contains $r$, or $\dist_G(P,Q)\le 2$.
\end{itemize}
This will prove \ref{mainthm}. (To see this, replace $S,T$ by $S\cup \{r\}, T\cup \{r\}$.)
We call $r$ the {\em root}. 

The construction for an $(\ell,1)$-block is easy: let $G'$ be a path of length $2\ell+1$ with ends $s,t$, and let $S=\{s\}$ and 
$T=\{t\}$; let $r$ be a new vertex, and let $G$ be obtained from $G'$ by adding $r$ as a vertex of degree zero. Then $(G,r,S,T)$
is an $(\ell,1)$-block.

\begin{figure}[h!]
\centering

\begin{tikzpicture}[scale=.5,auto=left]

\tikzstyle{every node}=[inner sep=1.5pt, fill=black,circle,draw]
\node (v2) at (2,0) {};
\node (v3) at (3,0) {};
\node (v4) at (4,0) {};
\node (v5) at (5,0) {};
\node (v6) at (6,0) {};
\node (v7) at (7,0) {};
\node (v8) at (8,0) {};
\node (v9) at (9,0) {};
\node (v10) at (10,0) {};
\node (v11) at (11,0) {};
\node (v12) at (12,0) {};
\node (v13) at (13,0) {};
\node (v14) at (14,0) {};
\node (v15) at (15,0) {};
\node (v16) at (16,0) {};
\node (v17) at (17,0) {};
\node (v18) at (18,0) {};
\node (v19) at (19,0) {};
\node (v20) at (20,0) {};
\node (v21) at (21,0) {};
\node (v22) at (22,0) {};
\node (v23) at (23,0) {};
\node (v24) at (24,0) {};
\node (v25) at (25,0) {};
\node (v26) at (26,0) {};
\node (v27) at (27,0) {};
\node (v28) at (28,0) {};
\node (v29) at (29,0) {};
\node (v30) at (30,0) {};
\node (v31) at (31,0) {};
\node (v32) at (32,0) {};
\node (v33) at (33,0) {};

\draw[dotted, thick] (v2)--(v33);

\node (u2) at (2.5,1) {};
\node (u4) at (4.5,1) {};
\node (u6) at (6.5,1) {};
\node (u8) at (8.5,1) {};
\node (u10) at (10.5,1) {};
\node (u12) at (12.5,1) {};
\node (u14) at (14.5,1) {};
\node (u16) at (16.5,1) {};
\node (u18) at (18.5,1) {};
\node (u20) at (20.5,1) {};
\node (u22) at (22.5,1) {};
\node (u24) at (24.5,1) {};
\node (u26) at (26.5,1) {};
\node (u28) at (28.5,1) {};
\node (u30) at (30.5,1) {};
\node (u32) at (32.5,1) {};

\draw[dotted,thick] (u2) to [bend right=20] (v2);
\draw[dotted,thick] (u4) to [bend right=20] (v3);
\draw[dotted,thick] (u6) to [bend right=20] (v5);
\draw[dotted,thick] (u8) to [bend right=20] (v7);
\draw[dotted,thick] (u10) to [bend right=20] (v9);
\draw[dotted,thick] (u12) to [bend right=20] (v11);
\draw[dotted,thick] (u14) to [bend right=20] (v13);
\draw[dotted,thick] (u16) to [bend right=20] (v15);
\draw[dotted,thick] (u18) to [bend right=20] (v17);
\draw[dotted,thick] (u20) to [bend right=20] (v19);
\draw[dotted,thick] (u22) to [bend right=20] (v21);
\draw[dotted,thick] (u24) to [bend right=20] (v23);
\draw[dotted,thick] (u26) to [bend right=20] (v25);
\draw[dotted,thick] (u28) to [bend right=20] (v27);
\draw[dotted,thick] (u30) to [bend right=20] (v29);
\draw[dotted,thick] (u32) to [bend right=20] (v31);

\draw[dotted,thick] (u2) to [bend left=20] (v4);
\draw[dotted,thick] (u4) to [bend left=20] (v6);
\draw[dotted,thick] (u6) to [bend left=20] (v8);
\draw[dotted,thick] (u8) to [bend left=20] (v10);
\draw[dotted,thick] (u10) to [bend left=20] (v12);
\draw[dotted,thick] (u12) to [bend left=20] (v14);
\draw[dotted,thick] (u14) to [bend left=20] (v16);
\draw[dotted,thick] (u16) to [bend left=20] (v18);
\draw[dotted,thick] (u18) to [bend left=20] (v20);
\draw[dotted,thick] (u20) to [bend left=20] (v22);
\draw[dotted,thick] (u22) to [bend left=20] (v24);
\draw[dotted,thick] (u24) to [bend left=20] (v26);
\draw[dotted,thick] (u26) to [bend left=20] (v28);
\draw[dotted,thick] (u28) to [bend left=20] (v30);
\draw[dotted,thick] (u30) to [bend left=20] (v32);
\draw[dotted,thick] (u32) to [bend left=20] (v33);

\node (t3) at (3.5,2) {};
\node (t7) at (7.5,2) {};
\node (t11) at (11.5,2) {};
\node (t15) at (15.5,2) {};
\node (t19) at (19.5,2) {};
\node (t23) at (23.5,2) {};
\node (t27) at (27.5,2) {};
\node (t31) at (31.5,2) {};

\draw (u2) -- (t3)--(u4);
\draw (u6) -- (t7)--(u8);
\draw (u10) -- (t11)--(u12);
\draw (u14) -- (t15)--(u16);
\draw (u18) -- (t19)--(u20);
\draw (u22) -- (t23)--(u24);
\draw (u26) -- (t27)--(u28);
\draw (u30) -- (t31)--(u32);

\node (s5) at (5.5,3) {};
\node (s13) at (13.5,3) {};
\node (s21) at (21.5,3) {};
\node (s29) at (29.5,3) {};

\draw (t3) -- (s5)--(t7);
\draw (t11) -- (s13)--(t15);
\draw (t19) -- (s21)--(t23);
\draw (t27) -- (s29)--(t31);

\node (r9) at (9.5,4) {};
\node (r25) at (25.5,4) {};

\draw (s5) -- (r9)--(s13);
\draw (s21) -- (r25)--(s29);

\node (q17) at (17.5,5) {};
\draw (r9) -- (q17)--(r25);

\tikzstyle{every node}=[]
\draw[above] (q17) node []           {$r$};
\draw[left] (u2) node []           {$S$};
\draw[left] (v2) node []           {$S$};
\draw[right] (v33) node []           {$T$};
\draw[right] (u32) node []           {$T$};

\end{tikzpicture}

\caption{The dotted curves represent paths of length $2\ell+1$.} \label{fig:counterexample}
\end{figure}
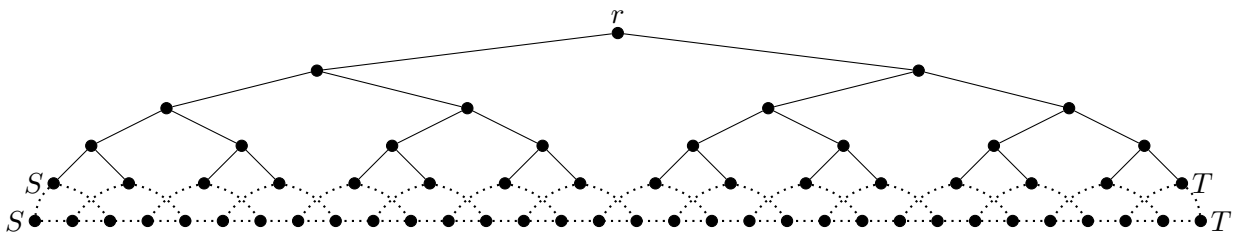

The construction for $(\ell,2)$-blocks was given in~\cite{counterex}, but here it is again (slightly modified for 
convenience), illustrated in figure~\ref{fig:counterexample}.
Each dotted curve in the figure represents a path of length $2\ell+1$, with interiors that are pairwise disjoint; let us call them ``dotted paths''. If we delete 
the interiors of the ``horizontal'' dotted paths at the bottom of the figure, 
we obtain a subdivision of a uniform binary tree, rooted at $r$, of depth 6 in the 
figure\footnote{The ``depth'' of a uniform binary tree is the number of vertices in paths from root to leaf.}; and to make an 
$(\ell,2)$-block, we need this tree to have depth at least $2\ell+2$. The path formed by the horizontal edges in the 
figure is called its {\em base path}; it starts in $S$ and ends in $T$, and has no other vertices in $S\cup T$. Let $s_0\in S$ and 
$t_0\in T$
be the vertices not in the base path. In the figure there are 32 vertices shown in the base path; they are the leaves of the binary 
tree. In general there are $2^{2\ell+1}$ such 
vertices, if the binary tree has depth $2\ell+2$. We call these the {\em anchors} of the base path. (The base path has many vertices 
that are not anchors, since it is a union of dotted paths). Let us number the anchors of the base path $v_1\CC v_n$ in order, where $v_1\in S$ and $v_n\in T$. 

An $(\ell,m+1)$-block is defined to be a certain combination of $(\ell,m)$-blocks, explained below.
Roughly, we replace each anchor of an $(\ell,2)$-block by a set of $m$ 
vertices, and replace each dotted path between anchors by an $(\ell,m)$-block, and then we identify all the roots.

\begin{figure}[h!]
\centering

\begin{tikzpicture}[scale=1,auto=left]

\tikzstyle{every node}=[inner sep=1.5pt, fill=black,circle,draw]
\node (v2) at (2,0) {};
\node (v3) at (3,0) {};
\node (v4) at (4,0) {};
\node (v5) at (5,0) {};
\node (v6) at (6,0) {};
\node (v7) at (7,0) {};
\node (v8) at (8,0) {};
\node (v9) at (9,0) {};
\node (v10) at (10,0) {};
\node (v11) at (11,0) {};
\node (v12) at (12,0) {};
\node (v13) at (13,0) {};
\node (v14) at (14,0) {};
\node (v15) at (15,0) {};
\node (v16) at (16,0) {};
\node (v17) at (17,0) {};

\node (w2) at (2,-1) {};
\node (w3) at (3,-1) {};
\node (w4) at (4,-1) {};
\node (w5) at (5,-1) {};
\node (w6) at (6,-1) {};
\node (w7) at (7,-1) {};
\node (w8) at (8,-1) {};
\node (w9) at (9,-1) {};
\node (w10) at (10,-1) {};
\node (w11) at (11,-1) {};
\node (w12) at (12,-1) {};
\node (w13) at (13,-1) {};
\node (w14) at (14,-1) {};
\node (w15) at (15,-1) {};
\node (w16) at (16,-1) {};
\node (w17) at (17,-1) {};

\node (u2) at (2.5,1) {};
\node (u4) at (4.5,1) {};
\node (u6) at (6.5,1) {};
\node (u8) at (8.5,1) {};
\node (u10) at (10.5,1) {};
\node (u12) at (12.5,1) {};
\node (u14) at (14.5,1) {};
\node (u16) at (16.5,1) {};

\draw[dotted,very thick] (u2) to [bend right=20] (v2);
\draw[dotted,very thick] (u4) to [bend right=20] (v3);
\draw[dotted,very thick] (u6) to [bend right=20] (v5);
\draw[dotted,very thick] (u8) to [bend right=20] (v7);
\draw[dotted,very thick] (u10) to [bend right=20] (v9);
\draw[dotted,very thick] (u12) to [bend right=20] (v11);
\draw[dotted,very thick] (u14) to [bend right=20] (v13);
\draw[dotted,very thick] (u16) to [bend right=20] (v15);

\draw[dotted,very thick] (u2) to [bend left=20] (v4);
\draw[dotted,very thick] (u4) to [bend left=20] (v6);
\draw[dotted,very thick] (u6) to [bend left=20] (v8);
\draw[dotted,very thick] (u8) to [bend left=20] (v10);
\draw[dotted,very thick] (u10) to [bend left=20] (v12);
\draw[dotted,very thick] (u12) to [bend left=20] (v14);
\draw[dotted,very thick] (u14) to [bend left=20] (v16);
\draw[dotted,very thick] (u16) to [bend left=20] (v17);

\node (t3) at (3.5,2) {};
\node (t7) at (7.5,2) {};
\node (t11) at (11.5,2) {};
\node (t15) at (15.5,2) {};

\draw (u2) -- (t3)--(u4);
\draw (u6) -- (t7)--(u8);
\draw (u10) -- (t11)--(u12);
\draw (u14) -- (t15)--(u16);

\node (s5) at (5.5,3) {};
\node (s13) at (13.5,3) {};

\draw (t3) -- (s5)--(t7);
\draw (t11) -- (s13)--(t15);

\node (r9) at (9.5,4) {};

\draw (s5) -- (r9)--(s13);

\tikzstyle{every node}=[]
\draw[above] (r9) node []           {$r_0$};
\draw[left] (u2) node []           {$S$};
\node[left] at (1.9,0) {$S$};

\draw[left] (w2) node []           {$S$};
\node[right] at (17.1,0) {$T$};

\draw[right] (u16) node []           {$T$};
\draw[right] (w17) node []           {$T$};
\draw[fill=gray!40] (v2) to [bend left=20] (w2) to [bend left=20] (w3) to [bend left=20] (v3) to [bend left=20] (v2);
\draw[fill=gray!40] (v3) to [bend left=20] (w3) to [bend left=20] (w4) to [bend left=20] (v4) to [bend left=20] (v3);
\draw[fill=gray!40] (v4) to [bend left=20] (w4) to [bend left=20] (w5) to [bend left=20] (v5) to [bend left=20] (v4);
\draw[fill=gray!40] (v5) to [bend left=20] (w5) to [bend left=20] (w6) to [bend left=20] (v6) to [bend left=20] (v5);
\draw[fill=gray!40] (v6) to [bend left=20] (w6) to [bend left=20] (w7) to [bend left=20] (v7) to [bend left=20] (v6);
\draw[fill=gray!40] (v7) to [bend left=20] (w7) to [bend left=20] (w8) to [bend left=20] (v8) to [bend left=20] (v7);
\draw[fill=gray!40] (v8) to [bend left=20] (w8) to [bend left=20] (w9) to [bend left=20] (v9) to [bend left=20] (v8);
\draw[fill=gray!40] (v9) to [bend left=20] (w9) to [bend left=20] (w10) to [bend left=20] (v10) to [bend left=20] (v9);
\draw[fill=gray!40] (v10) to [bend left=20] (w10) to [bend left=20] (w11) to [bend left=20] (v11) to [bend left=20] (v10);
\draw[fill=gray!40] (v11) to [bend left=20] (w11) to [bend left=20] (w12) to [bend left=20] (v12) to [bend left=20] (v11);
\draw[fill=gray!40] (v12) to [bend left=20] (w12) to [bend left=20] (w13) to [bend left=20] (v13) to [bend left=20] (v12);
\draw[fill=gray!40] (v13) to [bend left=20] (w13) to [bend left=20] (w14) to [bend left=20] (v14) to [bend left=20] (v13);
\draw[fill=gray!40] (v14) to [bend left=20] (w14) to [bend left=20] (w15) to [bend left=20] (v15) to [bend left=20] (v14);
\draw[fill=gray!40] (v15) to [bend left=20] (w15) to [bend left=20] (w16) to [bend left=20] (v16) to [bend left=20] (v15);
\draw[fill=gray!40] (v16) to [bend left=20] (w16) to [bend left=20] (w17) to [bend left=20] (v17) to [bend left=20] (v16);

\draw[dotted,very thick] (u2) to [bend right=40] (w2);
\draw[dotted,very thick] (u4) to [bend right=20] (w3);
\draw[dotted,very thick] (u6) to [bend right=20] (w5);
\draw[dotted,very thick] (u8) to [bend right=20] (w7);
\draw[dotted,very thick] (u10) to [bend right=20] (w9);
\draw[dotted,very thick] (u12) to [bend right=20] (w11);
\draw[dotted,very thick] (u14) to [bend right=20] (w13);
\draw[dotted,very thick] (u16) to [bend right=20] (w15);

\draw[dotted,very thick] (u2) to [bend left=20] (w4);
\draw[dotted,very thick] (u4) to [bend left=20] (w6);
\draw[dotted,very thick] (u6) to [bend left=20] (w8);
\draw[dotted,very thick] (u8) to [bend left=20] (w10);
\draw[dotted,very thick] (u10) to [bend left=20] (w12);
\draw[dotted,very thick] (u12) to [bend left=20] (w14);
\draw[dotted,very thick] (u14) to [bend left=20] (w16);
\draw[dotted,very thick] (u16) to [bend left=40] (w17);

\end{tikzpicture}

\caption{An $(\ell,3)$-block (except that the binary tree should have depth $2\ell+2$, and not just five as in the figure). The dotted
curves represent paths of length $2\ell+1$. Each gray area (together with $r_0$) is an $(\ell,2)$-block, and each contains neighbours of $r_0$ (not shown in the figure.)} \label{fig:weakcounterexample}

\end{figure}
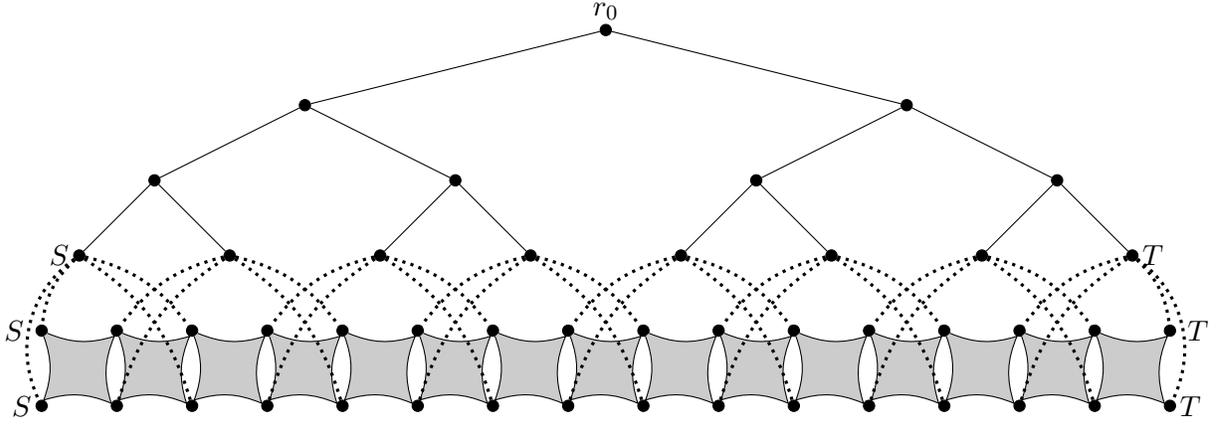

Here then, more exactly, is the inductive construction. We assume that we have an $(\ell,m)$-block $(H,r,S,T)$, and we will assemble copies of it
to make an $(\ell,m+1)$-block. Take an $(\ell,2)$-block $(H_0, r_0, S_0, T_0)$ as in Figure \ref{fig:counterexample}, 
and number the anchors of its base path
$v_1\LL v_n$, as described above. Let $S_0=\{v_1,s_0\}$ and $T_0=\{v_n,t_0\}$. For $1\le i\le n$, let $V_i$ be a set of $m$ new vertices; and for $1\le i<n$,
let $(H_i, r_0, V_i, V_{i+1})$ be a copy of $(H,r,S,T)$, where all vertices of
$H_i$ are new except those in $V_i\cup V_{i+1}\cup \{r_0\}$. Let $J$ be obtained from $H_0$ by deleting all vertices of its base path and all internal vertices of all its dotted 
paths (so $J$ is a
binary tree of depth $2\ell+1$). Let $G$ be the graph obtained from the union of $J, H_1\LL H_{n-1}$ by adding a path of length $2\ell+1$ between $u$ and each 
vertex of $V_i$, for each $u\in V(J)$ and each $v_i$ with $1\le i\le n$ such that $u,v_i$ are joined by a dotted path of $H_0$. 
We call these last {\em spines} of $G$.
Let $S=V_1\cup \{s_0\}$ and $T=V_n\cup \{y_0\}$. Then $(G,r_0,S,T)$ is an $(\ell, m+1)$-block.
(We illustrate this when $m=2$ in figure~\ref{fig:weakcounterexample}.)

It is easy to check that if $(G,r,S,T)$ is an $(\ell,m)$-block, then every two vertices in $\{r\}\cup S\cup T$ have distance at least $2\ell+1$, and now we will check the other properties mentioned in the bullets at the start of this section.
\begin{thm}\label{nodisjt}
If $(G,r_0,S,T)$ is an $(\ell,m+1)$-block, and $P,Q$ are paths between $S,T$, then either one of $P,Q$ contains $r_0$, or $\dist_G(P,Q)\le 2$.
\end{thm}
\Proof
We proceed by induction on $m$, and the result is true when $m=1$ by the result of \cite{counterex}, so we assume it is true for 
$(\ell,m)$-blocks. We use the notation given in the construction for an $(\ell,m+1)$-block.
Roughly, $G$ was obtained from $H_0$ by blowing up the vertices and edges of its base path; now we want to shrink them back to $H_0$, 
and carry $P,Q$ to some paths $P',Q'$ of $H_0$. More exactly, let $Z=V(J)\cup V_1\cupcup V_n$. The ends of $P$ both belong to $Z$;
so $P$ is the concatenation of a sequence of paths of $G$, each with distinct ends in $Z$ and no internal vertex in $Z$, say
$P_1\LL P_g$. Each $P_i$ is either an edge of $J$, or a spine, or a path of some $H_j$. We may assume that no $P_i$ contains $r_0$.
For $1\le i<n$, let $B_i$ be the dotted path of $H_0$ with ends $v_i, v_{i+1}$. 
For $1\le i\le g$, define a subgraph $P_i'$ of $H_0$ as follows.
\begin{itemize}
\item If $P_i$ is an edge of $J$ let $P_i'=P_i$. 
\item If $P_i$ is a spine with ends $u\in V(J)$ and some vertex in some $V_j$, let $P_i'$ be the dotted path
of $H_0$ between $u,v_i$. 
\item If $P_i$ is a path of some $H_j$ with both ends in $V_{j}$, let $P_i'$ be the one-vertex path with vertex $v_{j}$.
\item If $P_i$ is a path of some $H_j$ with both ends in $V_{j+1}$, let $P_i'$ be the one-vertex path with vertex $v_{j+1}$.
\item If $P_i$ is a path of some $H_j$ with one end in $V_j$ and the other in $V_{j+1}$, let $P_i'=B_j$.
\end{itemize}
It is not necessarily true that $P_1'\cupcup P_g'$ is a path of $H_0$ from $S_0$ to $T_0$, because it might pass through the same 
vertex more than once, and indeed $P_1'\LL P_g'$ need not all be distinct. But concatenating $P_1'\LL P_g'$ yields a walk from $S_0$ 
to $T_0$, and therefore $P_1'\cupcup P_g'$ includes a path $P'$ from $S_0$ to $T_0$. Define $Q'$ similarly, with $P$ replaced by $Q$. 
We may assume that no internal vertex of $P'$ or of $Q'$ belongs to $S_0\cup T_0$.
\\
\\
(1) {\em We may assume that there exists $i\in \{1\LL n\}$ such that $v_i\in V(P')\cap V(Q')$.}
\\
\\
Since $(H_0,r_0,S_0,T_0)$ is an $(\ell,2)$-block, and neither of $P',Q'$ contains $r_0$, it follows that 
$\dist_{H_0}(P',Q')\le 2$. 
Now there are two cases, depending whether $P',Q'$ are vertex-disjoint or not. If $P',Q'$ are vertex-disjoint, then since $\dist_{H_0}(P',Q')\le 2$, it follows that 
the corresponding path of $H_0$ joining them uses no dotted paths, and so is a path of $J$; and hence $\dist_G(P,Q)\le 2$, and the result holds.
So we may assume that there is some vertex $x\in V(P')\cap V(Q')$. If $x\in V(J)$, then $x\in V(P)\cap V(Q)$ and the result holds, so we 
may assume that $x\notin V(J)$. It remains to show that we may choose $x$ to be an anchor of $H_0$. If not, then $x$ belongs to 
the interior 
of some dotted path of $H_0$; but then the whole of that path belongs to 
$P'\cap Q'$, and one end of the path is an anchor. This proves (1).
\\
\\
(2) {\em We may assume that there exists $i\in \{1\LL n-1\}$ such that $B_i\subseteq P'\cap Q'$.}
\\
\\
By (1) we may assume that $v_i\in V(P'\cap Q')$ for some $i\in \{1\LL n\}$. 
Suppose first that $i = 1$.
Thus $v_1$ is the first vertex of both $P',Q'$. Neither of $P',Q'$ contain $s_0$, since 
 no internal vertex of $P'$ or of $Q'$ belongs to $S_0\cup T_0$; so both $P',Q'$ contain $B_1$ and the claim holds. 
Similarly the claim holds if $i=n$, so we assume that $2\le i\le n-1$. Now both $P',Q'$ include two of the three dotted paths of $H_0$
incident with $v_i$, and so include one in common. If that one is a spine, then its end in $V(J)$ belongs to $P'\cap Q'$ and we are done as
before; and otherwise 
one of $B_i, B_{i-1}\in P'\cap Q'$. This proves (2).

\bigskip

Consequently both $P,Q$ have subpaths with one end in $V_i$, one end in $V_{i+1}$, and all internal vertices in $H_i$. Since 
$(H_i, r_0, V_i, V_{i+1})$ is an $(\ell, m)$-block, and neither of $P,Q$ contains $r_0$, it follows that $\dist_{H_i}(P,Q)\le 2$,
and hence $\dist_G(P,Q)\le 2$. This proves \ref{nodisjt}.~\bbox

\begin{figure}[h!]
\centering

\begin{tikzpicture}[scale=2,auto=left]

\tikzstyle{every node}=[inner sep=1.5pt, fill=black,circle,draw]
\node (v2) at (2,0) {};
\node (v3) at (3,0) {};
\node (v4) at (4,0) {};
\node (v5) at (5,0) {};
\node (v6) at (6,0) {};
\node (v7) at (7,0) {};
\node (v8) at (8,0) {};
\node (v9) at (9,0) {};

\node (w2) at (2,-1) {};
\node (w3) at (3,-1) {};
\node (w4) at (4,-1) {};
\node (w5) at (5,-1) {};
\node (w6) at (6,-1) {};
\node (w7) at (7,-1) {};
\node (w8) at (8,-1) {};
\node (w9) at (9,-1) {};

\node (x2) at (2,-1/3) {};
\node (x3) at (3,-1/3) {};
\node (x4) at (4,-1/3) {};
\node (x5) at (5,-1/3) {};
\node (x6) at (6,-1/3) {};
\node (x7) at (7,-1/3) {};
\node (x8) at (8,-1/3) {};
\node (x9) at (9,-1/3) {};

\node (y2) at (2,-2/3) {};
\node (y3) at (3,-2/3) {};
\node (y4) at (4,-2/3) {};
\node (y5) at (5,-2/3) {};
\node (y6) at (6,-2/3) {};
\node (y7) at (7,-2/3) {};
\node (y8) at (8,-2/3) {};
\node (y9) at (9,-2/3) {};

\node (u2) at (2.5,1) {};
\node (u4) at (4.5,1) {};
\node (u6) at (6.5,1) {};
\node (u8) at (8.5,1) {};

\draw[dotted,very thick] (u4) to [bend right=20] (v3);
\draw[dotted,very thick] (u6) to [bend right=20] (v5);
\draw[dotted,very thick] (u8) to [bend right=20] (v7);

\draw[dotted,very thick] (u2) to [bend left=20] (v4);
\draw[dotted,very thick] (u4) to [bend left=20] (v6);
\draw[dotted,very thick] (u6) to [bend left=20] (v8);

\tikzstyle{every node}=[]

\draw[fill=gray!40] (v2) to [bend left=20] (x2) to [bend left=20] (y2) to [bend left=20] (w2) to [bend left=20] (w3) to 
[bend left=20] (y3) to [bend left=20] (x3) to [bend left=20] (v3) to [bend left=20] (v2);
\draw[fill=gray!40] (v3) to [bend left=20] (x3) to [bend left=20] (y3) to [bend left=20] (w3) to [bend left=20] (w4) to 
[bend left=20] (y4) to [bend left=20] (x4) to [bend left=20] (v4) to [bend left=20] (v3);
\draw[fill=gray!40] (v4) to [bend left=20] (x4) to [bend left=20] (y4) to [bend left=20] (w4) to [bend left=20] (w5) to 
[bend left=20] (y5) to [bend left=20] (x5) to [bend left=20] (v5) to [bend left=20] (v4);
\draw[fill=gray!40] (v5) to [bend left=20] (x5) to [bend left=20] (y5) to [bend left=20] (w5) to [bend left=20] (w6) to 
[bend left=20] (y6) to [bend left=20] (x6) to [bend left=20] (v6) to [bend left=20] (v5);
\draw[fill=gray!40] (v6) to [bend left=20] (x6) to [bend left=20] (y6) to [bend left=20] (w6) to [bend left=20] (w7) to 
[bend left=20] (y7) to [bend left=20] (x7) to [bend left=20] (v7) to [bend left=20] (v6);
\draw[fill=gray!40] (v7) to [bend left=20] (x7) to [bend left=20] (y7) to [bend left=20] (w7) to [bend left=20] (w8) to 
[bend left=20] (y8) to [bend left=20] (x8) to [bend left=20] (v8) to [bend left=20] (v7);
\draw[fill=gray!40] (v8) to [bend left=20] (x8) to [bend left=20] (y8) to [bend left=20] (w8) to [bend left=20] (w9) to 
[bend left=20] (y9) to [bend left=20] (x9) to [bend left=20] (v9) to [bend left=20] (v8);

\draw[dotted,very thick] (u4) to [bend right=20] (w3);
\draw[dotted,very thick] (u6) to [bend right=20] (w5);
\draw[dotted,very thick] (u8) to [bend right=20] (w7);

\draw[dotted,very thick] (u4) to [bend right=20] (x3);
\draw[dotted,very thick] (u6) to [bend right=20] (x5);
\draw[dotted,very thick] (u8) to [bend right=20] (x7);

\draw[dotted,very thick] (u4) to [bend right=20] (y3);
\draw[dotted,very thick] (u6) to [bend right=20] (y5);
\draw[dotted,very thick] (u8) to [bend right=20] (y7);

\draw[dotted,very thick] (u2) to [bend left=20] (w4);
\draw[dotted,very thick] (u4) to [bend left=20] (w6);
\draw[dotted,very thick] (u6) to [bend left=20] (w8);

\draw[dotted,very thick] (u2) to [bend left=20] (x4);
\draw[dotted,very thick] (u4) to [bend left=20] (x6);
\draw[dotted,very thick] (u6) to [bend left=20] (x8);

\draw[dotted,very thick] (u2) to [bend left=20] (y4);
\draw[dotted,very thick] (u4) to [bend left=20] (y6);
\draw[dotted,very thick] (u6) to [bend left=20] (y8);

\tikzstyle{every node}=[]
\node at (4,-1.25) {$V_{i-1}$};
\node at (5,-1.25) {$V_i$};
\node at (6,-1.25) {$V_{i+1}$};
\node at (7,-1.25) {$V_{i+2}$};

\node at (4.5,-.5) {$H_{i-1}$};
\node at (5.5,-.5) {$H_i$};
\node at (6.5,-.5) {$H_{i+1}$};

\draw[above] (u4) node []           {$u_{i-1}$};
\draw[above] (u6) node []           {$u_{i+1}$};

\draw[dashed,thick] (1.3,-1/2) -- (1.8,-1/2);
\draw[dashed,thick] (9.2,-1/2) -- (9.7,-1/2);

\end{tikzpicture}

\caption{Part of the construction for an $(\ell,m+1)$-block. The dotted
curves represent paths of length $2\ell+1$. Each gray area (together with $r_0$) is an $(\ell,m)$-block, and really has $m$ (not four as in the figure) vertices at either end. Each contains neighbours of $r_0$ (not shown in the figure.)} \label{fig:closeup}

\end{figure}
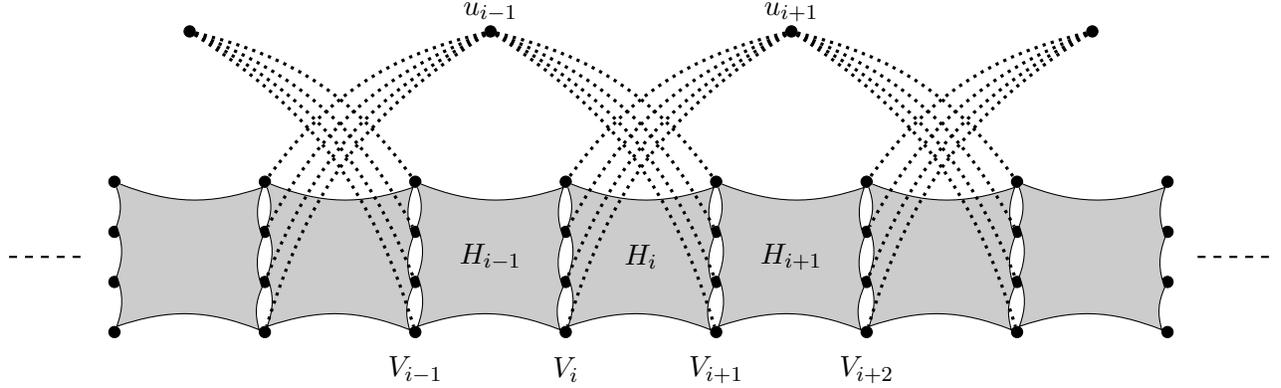

Next we need to show:

\begin{thm}\label{canthit}
If $(G,r_0,S,T)$ is an $(\ell,m+1)$-block,
and $X\subseteq V(G)$ with $|X|\le m$, then there is a path $P$ of $G$ between $S,T$
with $\dist_G(P,X\cup \{r_0\})>\ell$. 
\end{thm}
\Proof Again, we may assume the result holds for  $(\ell,m)$-blocks, and use the notation of the construction.
For each $x\in X$, let $C_x$ be the set of vertices $v$ of $G$ with $\dist_G(v,x)\le \ell$, 
let $C=\bigcup_{x\in X}C_x$, and let $C_0$ be the set of all $v$  with
$\dist_G(v,r_0)\le \ell$. We need to show that there is an $S-T$ path disjoint from $C\cup C_0$.
In fact we will prove that there is such a path within the union of the $H_i$'s and the spines. In Figure \ref{fig:closeup}
we show a section of the union of the $H_i$'s and spines. Let us label the vertices of $J$ that are ends of spines as in the figure.
(Thus, $s_0=u_1$, and $t_0=u_{n-1}$.)

We may assume that $m\ge 2$, and so every vertex in $S\cup T$ has degree at least three in $G$. We begin with:
\\
\\
(1) {\em We may assume that for every path $R$ of $G$ of length at least $2\ell+1$ in which all internal vertices have degree two 
in $G$, no internal vertex of $R$ belongs to $X$.}
\\
\\
If $x\in X$ is an internal vertex of such a path $R$, then by extending $R$ as much as possible, we may assume that 
both ends of $R$ have degree at least three. Let $R$ have ends $a,b$. Since $R$ has length at least $2\ell+1$, not all of $V(R)$
belongs to $C_x$, and we may assume that some vertex of $R$ between $x,b$ is not in $C_x$. Let $X'=(X\setminus \{x\})\cup \{a\}$.
Every path of $G$ within distance $\ell$ of $x$, and with no end in the interior of $R$, is within distance $\ell$ of $a$.
So if the result is true for $X'$ 
then it is true for $X$, so by repeating this at most $|X|$ times, we may assume there is no such $x$. This proves (1).
\\
\\
(2) {\em We may assume that for some $i\in \{1\LL n-1\}$, every path of $H_i\setminus r_0$ between $V_i, V_{i+1}$ intersects $C$.}
\\
\\
Suppose not, and for $1\le i\le n-1$ let $P_i$ be a path of $H_i\setminus r_0$ between $V_i, V_{i+1}$ that is vertex-disjoint
from $C$. For $2\le i<n$, the end of $P_{i-1}$ in $V_i$ and the end of $P_i$ in $V_i$ are joined by a path ($Q_i$ say) of $H_i$
of length $2\ell+1$ (or zero, if the two ends are equal) and all its internal vertices have degree two, so by (1), 
no vertex of $Q_i$ belongs to $X$.
Since its ends are not in $C$, no vertex of $Q_i$ is in $C$; and 
$$P_1\cup Q_2\cup P_2\cup Q_2\cupcup Q_{n-1}\cup P_{n-1}$$
is a path satisfying the theorem. This proves (2).

\bigskip

Let $i$ be as in (2). For each $x\in X$, $C_x$ contains at most one vertex in $V_i\cup V_{i+1}$, since they pairwise have distance 
at least $2\ell+1$. Consequently $C\cap V(H_i)$ is the union of at most $|X|$ balls of $H_i$ of radius $\ell$. But from the inductive hypothesis, 
for any set of at most $m-1$ such balls, there is a path of $H_i\setminus r_0$ between $V_i, V_{i+1}$ avoiding their union: 
so $C_x\cap V(H_i)\ne \emptyset$
for each $x\in X$. It follows that $X\subseteq V(H_{i-1})\cup V(H_i)\cup V(H_{i+1})$ (defining $H_0, H_n$ to be null), since 
no vertex of $X$ is in a spine by (1). Moreover, $V_{i-1}\cap C=\emptyset$ (again, defining  $V_0, V_{n+1}=\emptyset$), 
since $C_x\cap V(H_i)\ne \emptyset$ for each $x\in X$;
and similarly $V_{i+2}\cap C=\emptyset$. If $i$ is odd, the vertex called $u_i$ 
(see Figure \ref{fig:closeup}) exists, and provides a way to avoid $C$, via two spines incident with $u_i$; so we assume that 
$i$ is even. Thus $H_{i-1}, H_{i+1}$ 
exist. We may assume that $C$ does not contain $V_i$ (because if it does, then $i-1$ also satisfies (2) and is odd). 
Since $C_x\cap V(H_i)\ne \emptyset$ for each $x\in X$, and each of the sets $C_x\cap V_i\;(x\in X)$ as size at most one,
it follows that $C\cap V(H_{i-1})$ is the union of at most $|C\cap V_i|\le |V_i|-1= m-1$ balls of radius $\ell$, and so from 
the inductive hypothesis, 
there is a path of $H_{i-1}\setminus \{r_0\}$ between $V_{i-2}, V_{i-1}$ that is disjoint from $C$. But then we can extend this path 
via two spines incident with $u_{i+1}$ to provide a route that avoids $C$. This proves \ref{canthit}.~\bbox

From \ref{nodisjt} and \ref{canthit}, it follows that $(\ell,m-1)$-blocks provide a counterexample to \ref{weakcoarseconj} for each value of $\ell,m\ge 1$, and so this proves \ref{mainthm}.

\section{What remains?}
In view of this counterexample, is there anything even weaker that might be true and provide some sort of coarse extension of Menger's theorem? There are several possibilities to consider:

\subsubsection*{Bounded degree}
Imposing an upper bound on the maximum degree looks promising at first, because Gartland,  Korhonen and Lokshtanov~\cite{gartland} 
and Hendrey, Norin, Steiner, and Turcotte~\cite{turcotte} proved it with $c=2$: they proved
\begin{thm}\label{bddegree}
For every integer $\Delta\ge 1$ there exists $C>0$ with the following property.
Let $G$ be a graph, let $k\ge 1$ be an integer, and let $S,T\subseteq V(G)$; then either
\begin{itemize}
\item there are $k$ paths between $S,T$, pairwise at distance at least two; or
\item there is a set $X\subseteq V(G)$
with $|X|< kC$ such that every path between $S,T$ contains a vertex of $X$.
\end{itemize}
\end{thm}

But when $c\ge 3$, we can nearly make a counterexample. In the definition of an $(\ell,m+1)$-block, there are currently vertices of large degree:
the vertices of $J$ that are the ends of spines, and the root, and the vertices in the sets $V_i$. It is easy to modify the 
construction to keep the vertices in $V_i$ of degree at most three:  attach a leaf to each vertex in $S\cup T$, and call these new 
vertices $S$ and $T$ instead. (When we chain them together to make an $(\ell,m+1)$-block, they become of degree three instead of 
degree one, but that is fine.) For a vertex $u\in V(J)$ that is an end of spines, it is (in general) currently incident with $2m$ 
spines, 
and one edge of $J$. But we can partition these spines into two groups in the natural way, and replace each group with a binary tree
with root $u$ and leaves the corresponding set $V_i$; so that problem goes away as well. Thus, the only problem is the root vertex.
In summary, we can make a counterexample to \ref{weakcoarseconj} in which only one vertex has degree more than three, but we don't see how to fix this last vertex.
\subsection*{Bounded pathwidth}
Our counterexamples contain arbitrarily  large binary trees, and so have unbounded pathwidth. We will prove in a later paper~\cite{boundedpw} that
not only \ref{weakcoarseconj} but also \ref{conj} is true in graphs of bounded pathwidth:
\begin{thm}\label{boundedpw}
For all $c,k,w\ge 1$, there exists $\ell$ such that if $G$ has pathwidth at most $w$, and $S,T\subseteq V(G)$, then either:
\begin{itemize}
\item there are $k$ paths between $S,T$, pairwise at distance at least $c$; or
\item there is a set $X\subseteq V(G)$
with $|X|\le k-1$ such that every path between $S,T$ contains a vertex with distance at most $\ell$ from some member of $X$.
\end{itemize}
\end{thm}
Our counterexample to \ref{weakcoarseconj} contains subdivisions of arbitrarily large complete graphs, and indeed contains them 
with arbitrarily large ``fatness''
(see~\cite{agelos} for definition). It is easy to show (and several people have pointed out) that \ref{weakcoarseconj} is 
true for graphs of bounded tree-width. Is it true for every proper minor-closed class?
\subsection*{Planar graphs}
That suggests that we study whether \ref{weakcoarseconj} is true for planar graphs, but in that case, it seems possible that 
\ref{conj} itself is true. Indeed, \ref{conj} might be true for graphs of bounded genus: our counterexample to \ref{conj} 
in Figure \ref{fig:counterexample} has unbounded genus, because it contains arbitrarily many vertex-disjoint nonplanar subgraphs. 

We have begun to work on trying to prove \ref{conj} for planar graphs, 
and have been able to prove it in the case when the graph is planar and all the vertices in $S\cup T$ are incident with the 
infinite region~\cite{disccase}. (Despite its appearance, that was difficult!)  We think we can also prove it in the cylinder case,
when the graph is planar and some two regions include $S\cup T$: and perhaps there is an approach to the general coarse Menger 
conjecture for planar graphs that is like that in \cite{GM6,GM7}. But this needs much further work.

\end{document}